\documentclass[11pt]{article}
\usepackage{moreverb}
\usepackage{natbib}
\usepackage{setspace}
\usepackage{amsthm}
\usepackage{amssymb}
\usepackage{amsmath}
\usepackage{bm}
\usepackage{amsbsy}
\usepackage{accents}
\usepackage{palatino}
\usepackage{color}
\theoremstyle{plain}
\newtheorem{theorem}{Theorem}
\newtheorem{definition}{Definition}
\def\hgap{\ \ \ }
\def\bA{\boldsymbol{A}}
\def\bB{\boldsymbol{B}}
\def\bC{\boldsymbol{C}}
\def\bD{\boldsymbol{D}}

\def\bH{\boldsymbol{H}}
\def\bO{\boldsymbol{O}}
\def\bx{\boldsymbol{x}}
\def\ba{\boldsymbol{a}}
\def\bb{\boldsymbol{b}}
\def\bc{\boldsymbol{c}}
\def\bd{\boldsymbol{d}}
\def\bI{\boldsymbol{I}}
\def\bX{\boldsymbol{X}}
\def\bQ{\boldsymbol{Q}}
\def\bR{\boldsymbol{R}}
\def\bM{\boldsymbol{M}}
\def\bh{\boldsymbol{h}}
\DeclareMathOperator*{\blockdiag}{blockdiag}
\DeclareMathOperator*{\stack}{stack}
\def\bzero{\boldsymbol{0}}

\newcommand{\A}[1]{\bA_{#1}}
\newcommand{\AT}[1]{\bA_{#1}^T}
\newcommand{\Ainv}[1]{\bA^{#1}}
\newcommand{\ATinv}[1]{\bA^{#1 T}}
\def\smalldot{\mbox{\fontsize{0.1mm}{0.5em}\selectfont{$\bullet$}}}
\newcommand{\B}[1]{\bB_{#1}}
\newcommand{\BT}[1]{\bB_{#1}^T}
\newcommand*{\dt}[1]{\accentset{\mbox{\smalldot}}{#1}}
\newcommand*{\ddt}[1]{\accentset{\mbox{\smalldot\smalldot}}{#1}}
\newcommand{\dB}[1]{\dt{\boldsymbol{B}}_{#1}}
\newcommand{\ddB}[1]{\ddt{\boldsymbol{B}}_{#1}}
\newcommand{\dBT}[1]{\dt{\boldsymbol{B}}_{#1}^T}
\newcommand{\ddBT}[1]{\ddt{\boldsymbol{B}}_{#1}^T}
\newcommand{\C}[1]{\bC_{#1}}
\newcommand{\CT}[1]{\bC_{#1}^T}

\newcommand{\D}[1]{\bD_{#1}}
\newcommand{\DT}[1]{\bD_{#1}^T}
\newcommand{\dD}[1]{\dt{\boldsymbol{D}}_{#1}}
\newcommand{\dDT}[1]{\dt{\boldsymbol{D}}_{#1}^T}

\def\bigtimes{{\LARGE\mbox{$\times$}}}

\definecolor{ruppertgreen}{rgb}{0,.7,.25}

\def\TNbomega{\bc}
\def\TNbOmega{\bC}
\def\iCOMMAj{\,i,\,j}
\def\iCOMMAone{\,i,\,1}
\def\iCOMMAtwo{\,i,\,2}
\def\iCOMMAni{\,i,\,n_i}
\def\tcrORb#1{\textcolor{black}{#1}}
\newcommand\BibTeX{{\rmfamily B\kern-.05em \textsc{i\kern-.025em b}\kern-.08em
T\kern-.1667em\lower.7ex\hbox{E}\kern-.125emX}}
\begin{document}

\thispagestyle{empty}

\centerline{\Large\bf Solutions to Multilevel Sparse Matrix Problems}
\vskip7mm
\centerline{\normalsize\sc By Tui H. Nolan and Matt P. Wand}
\vskip5mm
\centerline{\textit{University of Technology Sydney}}
\vskip6mm
\centerline{11th March, 2020}

\vskip6mm
\centerline{\large\bf Abstract}
\vskip2mm

We define and solve classes of sparse matrix problems that arise in 
multilevel modeling and data analysis. The classes are indexed by
the number of nested units, with two-level problems corresponding
to the common situation in which data on level 1 units are grouped
within a two-level structure. We provide full solutions for two-level 
and three-level problems and their derivations provide blueprints 
for the challenging, albeit rarer in applications, 
higher level versions of the problem. Whilst our linear 
system solutions are a concise recasting of existing results, our
matrix inverse sub-block results are novel and facilitate
streamlined computation of standard errors in frequentist
inference as well as allowing streamlined mean field variational Bayesian 
inference for models containing higher level random effects.


\section{Introduction}

Higher level sparse matrices arise in statistical 
models for multilevel data, such as units
grouped according to geographical sub-regions or repeated
measures on medical study patients \citep[e.g.][]{Goldstein10}.
Other areas of statistics and econometrics that use
essentially the same types of models are longitudinal
data analysis \citep[e.g.][]{Fitzmaurice08}, panel data analysis \citep[e.g.][]{Baltagi13}
and small area estimation \citep[e.g.][]{Rao15}. Linear mixed models
\citep[e.g.][]{McCulloch08} are the main vehicle for modeling, fitting 
and inference. While they can be extended to generalised linear mixed
models to cater for skewed, categorical and count response data,
ordinary linear mixed models for Gaussian responses have the most
relevance for the sparse matrix results presented here.

Both frequentist and Bayesian estimates of the fixed and random 
effects can be expressed succinctly in terms of ridge regression-type
expressions involving design matrices \citep[e.g.][]{Henderson75}. 
However, typically the design matrices are sparse and 
na\"{\i}ve computation of the fixed and random effects 
estimates for large numbers of groups is inefficient and storage-greedy. 
For example, a random intercept linear mixed model for data 
with $1,000$ groups and $100$ observations per
group involves a random effects design matrix containing 100 million
entries of which 99.9\% are zeroes. Streamlined computation of the 
best linear unbiased predictors of the fixed and random effects
are well-documented with Section 2.2 of \citet{Pinheiro00} being
a prime example. The matrix algebraic notion of \emph{QR decomposition}
plays a central role in numerically stable least squares-based fitting 
of linear models \citep[e.g.][]{Gentle07} and also arises in the 
current context. It is important to note that such computations are 
performed \emph{after} estimates of the covariance matrix parameters 
have been obtained via approaches such as  minimum norm quadratic unbiased
estimation or restricted maximum likelihood. Streamlined computation
of covariance matrix estimates is tackled in, for example, 
\citet{Longford87}. Given the covariance matrix estimates,
implementation-ready matrix algebraic results for streamlined standard 
error calculations are not, to the best of our knowledge, present 
in the existing literature. These rely on efficient extraction of sub-blocks
of the inverses of potentially very large sparse symmetric matrices.
Presentation of these results, in the form of four theorems,
is our main novel contribution. In the interests of conciseness and
digestibility, we do not delve into the linear mixed model ramifications
here -- which are long-winded due to the various cases that require
separate treatment. This article is purely concerned with generic matrix
algebraic facts and, whilst motivated by statistical analysis, is
totally free of statistical concepts in its main results and derivations.
Ramifications for statistical inference are described in \citet{Nolan20}.

There is a related, but essentially non-overlapping, literature concerning
inversion of so-called \emph{arrowhead} matrices, which are \tcrORb{invertible} matrices
that have all entries equal to zero except for those on the main
diagonal and in one row and one column. \citet{Holubowski15} explain 
that such matrices ``often appear in areas of applied science 
and engineering such as head-positioning systems of hard disk drives or kinematic chains
of industrial robots'' and then develop a fast method for inversion 
of a generalisation of arrowhead matrices known as \emph{block arrowhead}
matrices. Other recent contributions of this type are \citet{Saberi14}
and \citet{Stanimirovic19}. It is important to note that,
whilst block arrowhead matrices coincide with two-level sparse matrices, the
class of problems treated here and their motivating statistical applications
are different from the central goal of the arrowhead matrix inversion literature. In 
multilevel sparse matrix problems the full matrix inverse is \emph{not}
of interest but, instead, inverse matrix sub-blocks matching the non-sparse sections 
of the original matrix. Arrowhead matrix inversion methods provide much more
than is required and, therefore, are overly slow for large sparse two-level
matrix problems.  In addition, multilevel sparse matrices beyond
the two-level case do not have block arrowhead forms.
In summary, multilevel sparse matrix problems are, in essence, 
distinct from arrowhead matrix inversion problems.

Motivated by applications described elsewhere \citep[e.g.][]{Nolan20},
our main focus in this article is the provision of full, implementable 
results for both two-level and three-level sparse matrix problems
with statistically relevant matrix inverse sub-blocks. Both 
general situations and least squares form situations, with QR decomposition 
enhancement, are covered. We believe that it is best to treat each higher level 
case separately. If a future application would benefit from the solution to the
four-level version of sparse matrix problems treated here then, whilst
notationally and algebraically challenging, our two-level and three-level 
derivations point the way to a solution.

We cover two-level sparse matrix problems in Section \ref{sec:2lev} 
and three-level sparse matrix problems in Section \ref{sec:3lev}.
Some concluding remarks are made in Section \ref{sec:conclusion}.

\section{Two-level Sparse Matrix Problems}\label{sec:2lev}

We begin by defining a two-level sparse matrix problem:

\begin{definition}
Let $\bA$ by a symmetric and invertible matrix of the form:
\begin{equation}
   \bA =
      \left[\arraycolsep=2.2pt\def\arraystretch{1.6}
         \begin{array}{c | c | c | c | c}
         \setstretch{4.5}
         \A{11} & \A{12,1} & \A{12,2} &\ \ \dots\ \ & \A{12,m} \\
         \hline
         \A{12,1}^T & \A{22,1} & \bO & \dots & \bO \\
         \hline
         \A{12,2}^T & \bO & \A{22,2} & \dots & \bO \\
         \hline
         \vdots & \vdots & \vdots & \ddots & \vdots \\
         \hline
         \A{12,m}^T & \bO & \bO & \dots & \A{22,m}
      \end{array} \right]
\label{2levMat}
\end{equation}
where the dimensions of the sub-blocks of $\bA$ are as follows:
$$\A{11}\ \text{is}\ p \times p\ \text{and, for each $1\le i\le m$,}\ \A{12,i}\ \text{is}\ p \times q
\ \text{and}\ \A{22,i}\ \text{is}\ q\times q.
$$
The two-level sparse matrix problem is defined to be:
\begin{enumerate}
\item[(I)] solve the linear system 
$$
\bA \bx = \ba,
$$
\item[(II)] obtain the sub-blocks of $\bA^{-1}$ corresponding to the  positions of the sub-blocks 
$\A{11}$, $\A{12,i}$ and $\A{22,i}$, $1\le i\le m$, in $\bA$.
\end{enumerate}
\label{def:twoLevDefn}
\end{definition}

If $p$ and $q$ are small relative to $m$ then matrices defined
by $\bA$ are sparse since, as $m\to\infty$, the leading term
of the fraction of non-zero entries of $\bA$ is $\{1+2(p/q)\}m^{-1}$.
In motivating statistical applications,
$m$ corresponds to sample size, which often is large, whilst $p$ and $q$ 
correspond to dimensions of model parameter spaces, which are moderate
in size. In the area of longitudinal data analysis typical values are 
$p=q=2$. The dimension variable $m$ matches the number of subjects 
in the longitudinal study, which may be in the hundreds or even 
thousands. If $m=1000$ then, for $p=q=2$, only about 0.3\% of the entries
of $\bA$ are non-zero. Throughout this article we assume that $p$ and $q$
are small relative to $m$. Similar remarks apply to the three-level
case treated in Section \ref{sec:3lev}.

Our solution to the two-level sparse matrix problem benefits from
the following notation for the sub-matrices of $\bA^{-1}$, $\ba$ and $\bx$:
\begin{equation}
   \bA^{-1} =
      \left[ 
         \arraycolsep=2.2pt\def\arraystretch{1.6}
         \begin{array}{c | c | c | c | c}
         \setstretch{4.5}
         \Ainv{11} & \Ainv{12,1} & \Ainv{12,2} &\ \ \dots\ \ & \Ainv{12,m} \\
         \hline
         \ATinv{12,1} & \Ainv{22,1} & \bigtimes & \dots & \bigtimes \\
         \hline
         \ATinv{12,2} & \bigtimes & \Ainv{22,2} & \dots & \bigtimes \\
         \hline
         \vdots & \vdots & \vdots & \ddots & \vdots \\
         \hline
         \ATinv{12,m} & \bigtimes & \bigtimes & \dots & \Ainv{22,m}
      \end{array} \right],
\qquad
   \ba \equiv
      \left[ 
       \arraycolsep=2.2pt\def\arraystretch{1.6}
         \begin{array}{c}
         \setstretch{4.5}
         \ba_{1} \\
         \hline
         \ba_{2,1} \\
         \hline
         \ba_{2,2} \\
         \hline
         \vdots \\
         \hline
         \ba_{2,m}
      \end{array} \right]
\qquad\text{and}\qquad
   \bx \equiv
      \left[\arraycolsep=2.2pt\def\arraystretch{1.6}
         \begin{array}{c}
         \setstretch{4.5}
         \bx_{1} \\
         \hline
         \bx_{2,1} \\
         \hline
         \bx_{2,2} \\
         \hline
         \vdots \\
         \hline
         \bx_{2,m}
      \end{array} \right]
\label{eq:dayFourTest}
\end{equation}
where $\bigtimes$ generically denotes sub-blocks of $\bA^{-1}$ which are in the same positions 
as the $\bO$ blocks in $\bA$. The dimensions of the sub-vectors of $\ba$ and $\bx$ are:
$$\mbox{both}\ \ba_1\ \mbox{and}\ \bx_1\ \mbox{are}\ p\times1\ \mbox{and, for $1\le i\le m$,}\ 
\mbox{both}\ \ba_{2,i}\ \mbox{and}\ \bx_{2,i}\ \mbox{are}\ q\times 1.
$$
Armed with the notation in (\ref{eq:dayFourTest}), part (II) of Definition \ref{def:twoLevDefn}
can be expressed as:
$$\mbox{obtain the matrix}\ \Ainv{11}\ \mbox{and, for each $1\le i\le m$, obtain the matrices}\ 
\Ainv{12,i}\ \mbox{and}\ \Ainv{22,i}.
$$

In applications involving multilevel data analysis the sub-blocks corresponding to the $\bigtimes$
symbols are usually not of interest since they correspond to between-group
covariances. On the other hand, the sub-blocks of $\bA^{-1}$ which are in the
same position as the non-zero sub-blocks of $\bA$ are required for obtaining
standard errors of within-group fits. In the case of mean field variational
inference, these sub-blocks are sufficient for both coordinate ascent 
and message passing optimal parameter computation  with minimal product 
restrictions. Details of how these sub-blocks of $\bA^{-1}$ are used in
linear mixed model inference are given in \citet{Nolan20}.

Theorem 1 provides a streamlined solution to this problem such that the 
number of operations is linear in $m$. An analogous expression for $|\bA|$, 
the determinant of $\bA$, is also provided.

\begin{theorem}
Consider the two-level sparse matrix problem given by Definition \ref{def:twoLevDefn}
\tcrORb{\normalsize and suppose that all $\A{22,i}$, $1\le i\le m$, are invertible}.
The solution to part (II) of Definition \ref{def:twoLevDefn} is:
{\setlength\arraycolsep{3pt}
\begin{eqnarray*}
\Ainv{11} &=& \left( \A{11} - \sum_{i=1}^m \A{12,i}\,\A{22,i}^{-1}\,\AT{12,i} \right)^{-1}\\[1ex]
\ \ \text{and}\ \ 
\Ainv{12,i} &=& -(\A{22,i}^{-1}\,\A{12,i}^T\,\Ainv{11})^T,
\ \ \Ainv{22,i} = \A{22,i}^{-1}(\bI-\A{12,i}^T\,\bA^{{12,i}}),\quad \text{$1\le i\le m$.}
\end{eqnarray*}
}
The determinant of $\bA$ is 
$$|\bA|=\left|\left(\Ainv{11}\right)^{-1}\right|\,\prod_{i=1}^m\big|\A{22,i}\big|.
$$
The solution to part (I) of Definition \ref{def:twoLevDefn} is:
$$\bx_1 = \Ainv{11} \left( \ba_1 - \sum_{i=1}^m \A{12,i}\,\A{22,i}^{-1}\,\ba_{2,i} \right)
\quad\text{and}\quad \bx_{2,i} =\A{22,i}^{-1}(\ba_{2,i} - \A{12,i}^T\,\bx_1), \quad 1 \le i \le m.
$$
\label{thm:firstThm}
\end{theorem}

\noindent
A proof of Theorem \ref{thm:firstThm} is given in Appendix \ref{sec:Th1}.

\subsection{Least Squares Form and QR-decomposition Enhancement}

In statistical applications involving linear mixed models, it is common for $\bA$ to admit
a \emph{least squares form} that lends itself to a QR decomposition-based solution. QR
decompositions of rectangular matrices are a numerically preferred method for solving
least squares problems. A QR-decomposition of a rectangular $n \times
p$ ($n \ge p$) matrix $\bX$ involves representing $\bX$ as
$$
\bX = \bQ \begin{bmatrix} \bR \\ \bO \end{bmatrix},
$$
where $\bQ$ is an $n \times n$ orthogonal matrix and $\bR$ is a $p \times p$
upper-triangular matrix. Such a decomposition affords computational stability for least
squares problems.

Suppose that $\bx$ is chosen to minimise the least squares criterion
\begin{equation}
   \Vert\bb - \bB \bx\Vert^2 \equiv (\bb - \bB \bx)^T (\bb - \bB \bx)
\label{leastSquares}
\end{equation}
for matrices
\begin{equation}
   \bB \equiv
      \left[ 
         \arraycolsep=2.2pt\def\arraystretch{1.6}
         \begin{array}{c | c | c | c | c}
         \setstretch{4.5}
         \bB_1 & \dB{1} & \bO & \dots & \bO \\
         \hline
         \bB_2 & \bO & \dB{2} & \dots & \bO \\
         \hline
         \vdots & \vdots & \vdots & \ddots & \vdots \\
         \hline
         \bB_m & \bO & \bO & \dots & \dB{m}
      \end{array} \right]
   \quad \text{and} \quad
   \bb \equiv
      \left[\arraycolsep=2.2pt\def\arraystretch{1.6} 
         \begin{array}{c}
         \setstretch{4.5}
         \bb_1 \\
         \hline
         \bb_2 \\
         \hline
         \vdots \\
         \hline
         \bb_m
      \end{array} \right]
\label{bmats}
\end{equation}
with sub-matrices and sub-vectors having dimensions:
$$\bB_i\ \mbox{is}\ n_i \times p,\quad \dB{i}\ \mbox{is}\ n_i \times q
\quad\mbox{and}\quad\bb_i\ \mbox{is}\ n_i \times 1
\quad\mbox{for $1\le i\le m$}
$$
\tcrORb{\normalsize such that $\bB$ is full rank.}
Then it is easily verified that the $\bx$ that minimises \eqref{leastSquares} 
is the solution to 
\begin{equation}
\bA\bx=\ba\quad\mbox{where}\quad
\bA = \bB^T \bB \quad\mbox{and}\quad \ba = \bB^T \bb.
\label{amats}
\end{equation}
Moreover, $\bA$ is a two-level sparse matrix of the form given by (\ref{2levMat}).
The non-zero sub-blocks of $\bA$ and the sub-vectors of $\ba$ are
\begin{equation*}
   \A{11} = \sum_{i=1}^m \bB_i^T \bB_i, \quad \ba_1 = \sum_{i=1}^m \BT{i} \bb_i
\end{equation*}
and
\begin{equation*}
   \A{12,i} = \bB_i^T \dB{i}, \quad \A{22,i} = \dBT{i} \dB{i}, \quad \ba_{2,i} = \dBT{i} \bb_i, \quad 1 \le i \le m.
\end{equation*}
The form of the matrix $\bB$ in (\ref{bmats}) arises in statistical models containing both fixed
effects and random effects with two-level structure (e.g. Goldstein, 2010). Full details on this motivational
connection are given in  Nolan \textit{et al.} (2020).

Theorem 2 extends Theorem 1 by employing a QR decomposition approach for
the purpose of numerical stability. Here, and later, we use the following notation
for matrices $\bM_1,\ldots,\bM_d$ each having the same number of columns:
$$\stack\limits_{1 \le i \le d}(\bM_i)\equiv
\left[
\begin{array}{c}
\bM_1\\
\vdots\\
\bM_d
\end{array}
\right].
$$

\begin{theorem}
Suppose that $\bA$ and $\ba$ admit the forms defined by \eqref{bmats} and \eqref{amats}
\tcrORb{\normalsize where $\bB$ is full rank.} Then the 
two-level sparse matrix problem may be solved using the following QR decomposition-based approach:
\begin{enumerate}

\item For $i=1,\ldots,m$:
\begin{enumerate}
\item Decompose $\dB{i} = \bQ_i \begin{bmatrix} \bR_i \\ \bzero \end{bmatrix}$ such 
that $\bQ_i^{-1} = \bQ_i^T$ and $\bR_i$ is upper-triangular.
\item Then obtain
   \begin{align*}
      &\bc_{0i} \equiv \bQ_i^T \bb_i, \quad \bc_{1i} \equiv \text{first q rows of} \ \bc_{0i}, 
      \quad \bc_{2i} \equiv \text{remaining rows of} \ \bc_{0i}, \\
      &\bC_{0i} \equiv \bQ_i^T \bB_i, \quad \bC_{1i} \equiv \text{first $q$ rows of} \ \bC_{0i}
         \quad \text{and} \quad \bC_{2i} \equiv \text{remaining rows of} \ \bC_{0i}.
   \end{align*}
\end{enumerate}
\item Decompose $\stack\limits_{1 \le i \le m} (\C{2i}) = \bQ \begin{bmatrix} \bR \\ \bO \end{bmatrix}$ 
such that $\bQ^{-1} = \bQ^T$ and $\bR$ is upper-triangular and let

   \begin{equation*}
      \TNbomega \equiv \text{first $p$ rows of} \ \bQ^T \left\{ \stack_{1 \le i \le m} (\bc_{2i}) \right\}.
   \end{equation*}

\item The solutions are
$$\bx_1 = \bR^{-1} \TNbomega,\quad\Ainv{11} = \bR^{-1} \bR^{-T}$$
and, for $1\le i\le m$,
$$\bx_{2,i} = \bR_i^{-1} (\bc_{1i} - \bC_{1i} \bx_1),\quad\Ainv{12,i} = -\Ainv{11} (\bR_i^{-1} \bC_{1i})^T,
\quad \Ainv{22,i} = \bR_i^{-1}(\bR_i^{-T}-\bC_{1i}\Ainv{12,i}).
$$
\item The determinant of $\bA$ is
$$|\bA|=\left\{\big(\mbox{product of the diagonal entries of $\bR$}\big)\,\prod_{i=1}^m\, 
\big(\mbox{product of the diagonal entries of $\bR_i$}\big)\right\}^2.$$
\end{enumerate}
\label{thm:secondThm}
\end{theorem}

\noindent
A proof of Theorem \ref{thm:secondThm} is in Appendix \ref{sec:Th2}.

\noindent
\textsl{Remarks:}
\begin{enumerate}
\item In Theorem \ref{thm:secondThm}, Step 1 involves determination of $m$ upper triangular matrices 
$\bR_i$, $1 \le i \le m$, via QR-decomposition which is a standard procedure within most computing 
environments. Each of the matrix inversions in Step 3 involve $\bR_i^{-1}$, which can be achieved 
rapidly via back-solving.
\item Calculations such as $\bQ_i^T\bb_i$ do not require storage of $\bQ_i$ and ordinary 
matrix multiplication. Standard matrix algebraic programming languages are such that information
concerning $\bQ_i$ is stored in a compact form from which matrices such as $\bQ_i^T\bb_i$ 
can be efficiently obtained.
\item Pinheiro \& Bates (2000; Section 2.2) make use of this QR decomposition-based approach for 
fitting two-level linear mixed models. However, their descriptions are restricted to the $\bx_1$ 
and $\bx_{2,i}$ formulae, not those for the sub-blocks of $\bA^{-1}$.
\end{enumerate}

Table \ref{tab:timeStudyComp} summarises the results of a numerical study 
for which $100$ $\bA$ matrices of the form (\ref{2levMat}) were randomly generated
and the solution to the streamlined two-level sparse matrix problem using Theorem \ref{thm:secondThm}
was compared to  the na\"{\i}ve solution, where the sparse structure
is ignored, for increasingly large versions of the problem.
Throughout the study $p=q=2$ and $m$ ranged over the
set $\{100,200,400,800,1600\}$. The $n_i$ values were generated uniformly on the set 
$\{30,\ldots,60\}$. To aid maximal speed, both approaches were implemented in the 
low-level language \texttt{Fortran 77}, with \texttt{LINPACK} subroutines used
for numerical linear algebra. The na\"{\i}ve approach involved matrix inversion via 
Gaussian elimination. The study was run on a \texttt{MacBook Air} 
laptop computer with a 2.2 gigahertz processor and 8 gigabytes of random access memory.
The table lists the average and standard deviation times in seconds across the 
$100$ replications. 
%%%%%%%%%%%%
\vfill\eject
%%%%%%%%%%%%

%%%%%%%%%%%%%%%%%%%%%%%%%%%%%%%%%%%%%%%%%%%%%%%%%%%%%%%%%%%%%%%%%%%%%%%%%%%
\begin{table}[h]
\begin{center}
\begin{tabular}{cccc}
\hline\\[-2.5ex]
$m$                & na\"{\i}ve & streamlined  &  na\"{\i}ve/streamlined\\
\hline\\[-1.5ex]
100               &   0.103 (0.0173)&    0.00253 (0.000541)  &   40.8\\
200               &   0.768 (0.0214)&    0.00429 (0.000574)  &  179.0\\
400               &   6.050 (0.1110)&    0.00798 (0.000586)  &  758.0\\
800               &  48.300 (0.7830)&    0.01530 (0.001650)  & 3150.0\\
1600              & 386.000 (6.0200)&    0.03000 (0.001900)  &12900.0\\
\hline
\end{tabular}
\end{center}
\caption{\textit{Averages (standard deviations) of elapsed computing times in seconds for 
solving the two-level sparse matrix problem na\"{\i}vely versus with the streamlined
approach provided by Theorem \ref{thm:secondThm}. The fourth column lists the ratios of the median
computing times.}}
\label{tab:timeStudyComp}
\end{table}
%%%%%%%%%%%%%%%%%%%%%%%%%%%%%%%%%%%%%%%%%%%%%%%%%%%%%%%%%%%%%%%%%%%%%%%%%

From Table \ref{tab:timeStudyComp} we see that the streamlined solutions
are delivered in a small fraction of a second, even for very large versions
of the problem. By contrast, na\"{\i}ve computation takes several minutes
for the $m=1600$ case and is more than ten thousand times slower. 

\section{Three-level Sparse Matrix Problems}\label{sec:3lev}

Three-level sparse matrix problems are such that the $\bA$ matrix is symmetric and invertible
but the two-level structure is repeated down the main diagonal. The notation for the general
case becomes difficult to digest, so we start with a concrete example of such an $\bA$:
\begin{equation}
\bA =
\left[ \arraycolsep=2.2pt\def\arraystretch{1.6} 
   \begin{array}{c | c | c | c | c | c | c | c}
   \setstretch{4.5}
   \A{11} & \A{12,1} & \A{12,11} & \A{12,12} & \A{12,2} & \A{12,21} & \A{12,22} & \A{12,23} \\
   \hline
   \A{12,1}^T & \A{22,1} & \A{12,1,1} & \A{12,1,2} & \bO & \bO & \bO & \bO \\
   \hline
   \A{12,11}^T & \A{12,1,1}^T & \A{22,11} & \bO & \bO & \bO & \bO & \bO \\
   \hline
   \A{12,12}^T & \A{12,1,2}^T & \bO & \A{22,12} & \bO & \bO & \bO & \bO \\
   \hline
   \A{12,2}^T & \bO & \bO & \bO & \A{22,2} & \A{12,2,1} & \A{12,2,2} & \A{12,2,3} \\
   \hline
   \A{12,21}^T & \bO & \bO & \bO & \A{12,2,1}^T & \A{22,21} & \bO & \bO \\
   \hline
   \A{12,22}^T & \bO & \bO & \bO & \A{12,2,2}^T & \bO & \A{22,22} & \bO \\
   \hline
   \A{12,23}^T & \bO & \bO & \bO & \A{12,2,3}^T & \bO & \bO & \A{22,23}
\end{array} \right].
\label{eq:AthreeLevExamp}
\end{equation}
The general three-level sparse matrix problem is defined by:
\begin{definition}
Let $\bA$ by a symmetric and invertible matrix with partitioning
$$\bA=\left[\begin{array}{cc}\A{11}   & \A{12}\\[1ex]
                             \A{12}^T & \A{22}
\end{array}\right]
$$
where $\A{11}$ is $p\times p$ and
$$\bA_{12}=\left\{\left[\begin{array}{c| c|c|c}\A{12,i}& \A{12,ij} 
& \dots & \A{12,in_i} \end{array} \right] \right\}_{1 \le i \le m}
$$
such that for each $1 \le i \le m$, $\A{12,i}$ 
is $p \times q_1$, and for each $1 \le j \le n_i$, $\A{12,ij}$ is $p \times q_2$.
The lower right block is
$$\bA_{22}=\blockdiag_{1\le i\le m}
      \left(\left[\arraycolsep=2.2pt\def\arraystretch{1.6}
         \begin{array}{c | c | c | c | c}
         \setstretch{4.5}
         \A{22,i} & \A{12,\iCOMMAone} & \A{12,\iCOMMAtwo} &\ \ \dots\ \ & \A{12,\iCOMMAni} \\
         \hline
         \A{12,\iCOMMAone}^T & \A{22,i1} & \bO & \dots & \bO \\
         \hline
         \A{12,\iCOMMAtwo}^T & \bO & \A{22,i2} & \dots & \bO \\
         \hline
         \vdots & \vdots & \vdots & \ddots & \vdots \\
         \hline
         \A{12,\iCOMMAni}^T & \bO & \bO & \dots & \A{22,in_i}
      \end{array} \right]\right)
$$
where, for each $1 \le i \le m$, $\A{22,i}$ is $q_1 \times q_1$, and for each 
$1 \le j \le n_i$, $\A{12,\iCOMMAj}$ is $q_1 \times q_2$ and $\A{22,ij}$ is $q_2\times q_2$.
The three-level sparse matrix problem is defined to be: 
\begin{enumerate}
\item[(I)] solve the linear system 
$$
\bA \bx = \ba,
$$
\item[(II)] obtain the sub-blocks of $\bA^{-1}$ corresponding to the positions of the sub-blocks 
$\A{11}$, $\A{12,i}$, $\A{22,i}$, $1\le i\le m$, and $\A{12,ij},\ \A{12,\iCOMMAj},\ \A{22,ij}$, 
$1\le i\le m$, $1\le j\le n_i$, in $\bA$.
\end{enumerate}
\label{def:threeLevDefn}
\end{definition}

\noindent 
For the example three-level sparse matrix given by (\ref{eq:AthreeLevExamp}), 
we have $m=2$, $n_1=2$ and $n_2=3$.
To enhance digestibility we will use these values for $m$, $n_1$ and $n_2$ throughout 
our discussion regarding the three-level sparse matrix problem. However, this can be 
easily generalised to any three-level sparse matrix, where $m$ and $\{ n_i \}_{1 \le i \le m}$ 
are arbitrary.

The solution to the thee-level sparse matrix problem benefits from
the following notation for the sub-matrices of $\bA^{-1}$, $\ba$ and $\bx$:
\begin{equation}
\bA^{-1} =
\left[\arraycolsep=2.2pt\def\arraystretch{1.6}
   \begin{array}{c | c | c | c | c | c | c | c}
 \setstretch{4.5}
   \Ainv{11} & \Ainv{12,1} & \Ainv{12,11} & \Ainv{12,12} & \Ainv{12,2} & \Ainv{12,21} & \Ainv{12,22} & \Ainv{12,23} \\
   \hline
   \ATinv{12,1} & \Ainv{22,1} & \Ainv{12,1,1} & \Ainv{12,1,2} & \bigtimes & \bigtimes & \bigtimes & \bigtimes \\
   \hline
   \ATinv{12,11} & \ATinv{12,1,1} & \Ainv{22,11} & \bigtimes & \bigtimes & \bigtimes & \bigtimes & \bigtimes \\
   \hline
   \ATinv{12,12} & \ATinv{12,1,2} & \bigtimes & \Ainv{22,12} & \bigtimes & \bigtimes & \bigtimes & \bigtimes \\
   \hline
   \ATinv{12,2} & \bigtimes & \bigtimes & \bigtimes & \Ainv{22,2} & \Ainv{12,2,1} & \Ainv{12,2,2} & \Ainv{12,2,3} \\
   \hline
   \ATinv{12,21} & \bigtimes & \bigtimes & \bigtimes & \ATinv{12,2,1} & \Ainv{22,21} & \bigtimes & \bigtimes \\
   \hline
   \ATinv{12,22} & \bigtimes & \bigtimes & \bigtimes & \ATinv{12,2,2} & \bigtimes & \Ainv{22,22} & \bigtimes \\
   \hline
   \ATinv{12,23} & \bigtimes & \bigtimes & \bigtimes & \ATinv{12,2,3} & \bigtimes & \bigtimes & \Ainv{22,23}
\end{array} \right],
\label{eq:stevie}
\end{equation}
$$
\ba \equiv
      \left[ \arraycolsep=2.2pt\def\arraystretch{1.6}
         \begin{array}{c}
         \setstretch{4.5}
         \ba_{1} \\
         \hline
         \ba_{2,1} \\
         \hline
         \ba_{2,11} \\
         \hline
         \ba_{2,12} \\
         \hline
         \ba_{2,2} \\
         \hline
         \ba_{2,21} \\
         \hline
         \ba_{2,22} \\
         \hline
         \ba_{2,23}
      \end{array} \right] 
\quad\mbox{and}\quad 
   \bx \equiv
      \left[ 
         \arraycolsep=2.2pt\def\arraystretch{1.6}
         \begin{array}{c}
         \setstretch{4.5}
         \bx_{1} \\
         \hline
         \bx_{2,1} \\
         \hline
         \bx_{2,11} \\
         \hline
         \bx_{2,12} \\
         \hline
         \bx_{2,2} \\
         \hline
         \bx_{2,21} \\
         \hline
         \bx_{2,22} \\
         \hline
         \bx_{2,23}
      \end{array} \right].
$$
The dimensions of the partitioned vectors are:
\begin{itemize}
\item $\ba_1$ and $\bx_1$ are $p \times 1$ vectors;
\item for each $1 \le i \le m$, $\ba_{2,i}$ and $\bx_{2,i}$ are $q_1 \times 1$ vectors;
\item for each $1 \le i \le m$ and $1 \le j \le n_i$, $\ba_{2,ij}$ and $\bx_{2,ij}$ are
$q_2 \times 1$ vectors.
\end{itemize}

\noindent As in Section \ref{sec:2lev}, $\bigtimes$ denotes the blocks of $\bA^{-1}$ that are not of interest.
Using the notation exemplified by (\ref{eq:stevie}), part (II) of Definition \ref{def:threeLevDefn}
can be expressed as:
$$
\begin{array}{c}
\mbox{obtain the matrix}\ \Ainv{11},\ \mbox{for each $1\le i\le m$, obtain the matrices}\ 
\Ainv{12,i},\ \Ainv{22,i}\ \mbox{and}\\
\mbox{for each $1\le i\le m$, $1\le j\le n_i$ obtain the matrices}
\ \Ainv{12,ij},\ \Ainv{12,\iCOMMAj},\ \Ainv{22,ij}.
\end{array}
$$
Theorem 3 presents the solution for any matrix that has
the same sparsity structure as $\bA$  when $m$ and $\{ n_i \}_{1 \le i \le m}$ are arbitrary.
An analogous expression for $|\bA|$ is also provided.

\begin{theorem}
Consider the three-level sparse matrix problem given by Definition \ref{def:threeLevDefn}
\tcrORb{\normalsize and suppose that all $\A{22,ij}$, $1\le i\le m$, $1\le j\le n_i$,  are invertible}.
For $1\le i\le m$, define
{\setlength\arraycolsep{3pt}
\begin{eqnarray*}
\bh_{2,i}&\equiv&\ba_{2,i} - \sum_{j=1}^{n_i}\,\A{12,\iCOMMAj}\,\A{22,ij}^{-1}\,\ba_{2,ij},\quad
\bH_{12,i}\equiv \A{12,i}-\sum_{j=1}^{n_i}\A{12,ij}\,\A{22,ij}^{-1}\,\AT{12,\iCOMMAj}\\
\quad\text{and}\quad
\bH_{22,i}&\equiv&\A{22,i}-\sum_{j=1}^{n_i}\A{12,\iCOMMAj}\,\A{22,ij}^{-1}\,\AT{12,\iCOMMAj}
\end{eqnarray*}
}
\tcrORb{\normalsize and suppose that all $\bH_{22,i}$, $1\le i\le m$,  are invertible}.
The solution to part (II) of Definition \ref{def:threeLevDefn} is:
\begin{align*}
   \Ainv{11} &=
      \left(
         \A{11}- \sum_{i=1}^m \sum_{j=1}^{n_i}\A{12,ij}\,\A{22,ij}^{-1}\,\A{12,ij}^T
         - \sum_{i=1}^m \bH_{12,i}\,\bH_{22,i}^{-1}\,\bH_{12,i}^T         
      \right)^{-1},\\[1ex]
   \Ainv{12,i} &= -(\bH_{22,i}^{-1}\,\bH_{12,i}^T\,\Ainv{11})^T,\quad 
   \Ainv{22,i} = \bH_{22,i}^{-1}(\bI-\,\bH_{12,i}^T\bA^{12,i}),\quad \text{$1 \le i \le m$,}\\[1ex] 
   \Ainv{12,ij} &= -\big\{\A{22,ij}^{-1}\big(\A{12,ij}^T\,\Ainv{11}\ + \A{12,\iCOMMAj}^T\,\Ainv{12,i\,T}\big)\big\}^T,\\[1ex]
   \Ainv{12,\iCOMMAj}&= -\big\{\A{22,ij}^{-1}\big(\AT{12,ij}\,\Ainv{12,i}+\A{12,\iCOMMAj}^T\,\Ainv{22,i}\big)\big\}^T,\\[1ex]
\Ainv{22,ij}&= \A{22,ij}^{-1}\big(\bI-\bA_{12,ij}^T\,\bA^{12,ij}-\bA_{12,\iCOMMAj}^T\,\bA^{12,\iCOMMAj}\big),
\quad\text{$1 \le i \le m$,\ $1 \le j \le n_i$.}
\end{align*}
The determinant of $\bA$ is 
$$
|\bA|=
\left|\left(\Ainv{11}\right)^{-1}\right|
\,\prod_{i=1}^m\left(\left|\bA_{22,i}-\sum_{j=1}^{n_i}\bA_{12,\iCOMMAj}\,
\A{22,ij}^{-1}\,\bA_{12,\iCOMMAj}^T\right|
\prod_{j=1}^{n_i}\big|\A{22,ij}\big|\right).
$$
The solution to part (I) of Definition \ref{def:threeLevDefn} is:
\begin{align*}
   \bx_1 &=
      \Ainv{11}\left(
         \ba_1-\sum_{i=1}^m\bH_{12,i}\,\bH_{22,i}^{-1}\,\bh_{2,i}
              -\sum_{i=1}^m \sum_{j=1}^{n_i}\A{12,ij}\,\A{22,ij}^{-1}\,\ba_{2,ij}
              \right),\\[1ex]
   \bx_{2,i} &=
      \bH_{22,i}^{-1}\left(\bh_{2,i}- \bH_{12,i}^T\,\bx_1
      \right), \quad \text{$1 \le i \le m$}, \\[1ex]
   \bx_{2,ij} &= \
      \A{22,ij}^{-1}\left(\ba_{2,ij} - \A{12,ij}^T\,\bx_1 - \A{12,\iCOMMAj}^T\,\bx_{2,i}
      \right), \quad \text{$1\le i\le m$, $1\le j\le n_i$}.
\end{align*}
\label{thm:thirdThm}
\end{theorem}

\noindent
A proof of Theorem \ref{thm:thirdThm} is in Appendix \ref{sec:Th3}.

\subsection{Least Squares Form and QR-decomposition Enhancement}

The three-level sparse matrix problem also lends itself to QR-decomposition enhancement. 
For the special case of $\bA$, with $m=2$, $n_1=2$ and $n_2=3$, the least squares criterion 
has the form (\ref{leastSquares}) with
\begin{equation}
   \bB \equiv
      \left[\arraycolsep=2.2pt\def\arraystretch{1.6}
         \begin{array}{c | c | c | c | c | c | c | c}
         \setstretch{4.5}
         \B{11} & \dB{11} & \ddB{11} & \bO & \bO & \bO & \bO & \bO \\
         \hline
         \B{12} & \dB{12} & \bO & \ddB{12} & \bO & \bO & \bO & \bO \\
         \hline
         \B{21} & \bO & \bO & \bO & \dB{21} & \ddB{21} & \bO & \bO \\
         \hline
         \B{22} & \bO & \bO & \bO & \dB{22} & \bO & \ddB{22} & \bO \\
         \hline
         \B{23} & \bO & \bO & \bO & \dB{23} & \bO & \bO & \ddB{23}
      \end{array} \right]
   \quad \text{and} \quad
   \bb \equiv
      \left[\arraycolsep=2.2pt\def\arraystretch{1.6}
         \begin{array}{c}
         \setstretch{4.5}
         \bb_{11} \\
         \hline
         \bb_{12} \\
         \hline
         \bb_{21} \\
         \hline
         \bb_{22} \\
         \hline
         \bb_{23}
      \end{array} \right].
\label{Bmats3lev}
\end{equation}
For general values of $m$ and $\{n_i\}_{1\le i\le m}$, the forms of $\bB$ and  $\bb$ are
\begin{equation}
\bB\equiv\Big[\stack\limits_{1\le i\le m}\Big\{\stack\limits_{1\le j\le n_i}(\B{ij})\Big\}\ \Big\vert
\blockdiag\limits_{1\le i\le m}\Big\{\big[\stack\limits_{1\le j\le n_i}(\dB{ij})\ \big\vert
\  \blockdiag\limits_{1\le j\le n_i}(\ddB{ij})\big]\Big\} \Big]
\label{eq:threeLevB}
\end{equation}
and
\begin{equation}
\bb\equiv\stack\limits_{1\le i\le m}\Big\{\stack\limits_{1\le j\le n_i}(\bb_{ij})\Big\}.
\label{eq:threeLevSmallb}
\end{equation}
For each $1 \le i \le m$ and $1 \le j \le n_i$, the dimensions of the sub-blocks of $\bB$ and $\bb$ are:
$$\B{ij}\ \text{is}\ o_{ij}\times p,\quad\dB{ij}\ \text{is}\ o_{ij}\times q_1,
\quad\ddB{ij}\ \text{is}\ o_{ij}\ \times q_2
\quad\mbox{and}\quad\bb_{ij}\ \text{is}\ o_{ij}\times 1.
$$
\tcrORb{\normalsize We also assume that $\bB$ is full rank.}
Then the $\bx$ that minimises \eqref{leastSquares} is the solution to the three-level 
sparse linear system with
\begin{equation}
\bA = \bB^T \bB \quad \text{and} \quad \ba = \bB^T \bb.
\label{eq:AaDefnThreeLev}
\end{equation}
For general $m$ and $\{ n_i \}_{1 \le i \le m}$ the non-zero components of $\bA$ and the sub-vectors of $\ba$ are,
for $1\le i\le m$,
{\setlength\arraycolsep{3pt}
\begin{eqnarray*}
\A{11} &=& \sum_{i=1}^m \sum_{j=1}^{n_i} \BT{ij} \B{ij},\quad 
\ba_1 = \sum_{i=1}^m \sum_{j=1}^{n_i} \BT{ij} \bb_{ij},\quad
\A{22,i} = \sum_{j=1}^{n_i} \dBT{ij} \dB{ij},\\[1ex] 
\A{12,i} &=& \sum_{j=1}^{n_i} \BT{ij} \dB{ij},\quad
\ba_{2,i} = \sum_{j=1}^{n_i} \dBT{ij} \bb_{ij}
\end{eqnarray*}
}
and
$$
\A{22,ij} = \ddBT{ij} \ddB{ij}, \quad \A{12,\iCOMMAj} = \dBT{ij} \ddB{ij},\quad
\A{12,ij} = \BT{ij} \ddB{ij}, \quad \ba_{2,ij} = \ddBT{ij} \bb_{ij}, \quad 1 \le i \le m,\ 1 \le j \le n_i.
$$
The form of the matrix $\bB$ in (\ref{eq:threeLevB}) arises in statistical models containing both
fixed effects and random effects with three-level structure (e.g. Goldstein, 2010). Fuller details on this
connection are given in  Nolan \textit{et al.} (2020).

Theorem 4 provides a QR decomposition enhancement of Theorem 3 for the least squares forms situation.

\begin{theorem}
Suppose that $\bA$ and $\ba$ admit the least squares forms defined by 
(\ref{eq:threeLevB})--(\ref{eq:AaDefnThreeLev}) \tcrORb{\normalsize where $\bB$ is full rank}. 
Then the three-level sparse matrix problem may be solved using the following QR decomposition-based approach:
\begin{enumerate}
\item For $i=1,\ldots,m$:
\begin{enumerate}
\item For $j=1,\ldots,n_i$:
\begin{enumerate}
\item Decompose $\ddB{ij} = \bQ_{ij} \begin{bmatrix} \bR_{ij} \\ \bzero \end{bmatrix}$ 
such that $\bQ_{ij}^{-1} = \bQ_{ij}^T$ and $\bR_{ij}$ is upper-triangular.
\item Then obtain
\begin{align*}
      &\bd_{0ij} \equiv \bQ_{ij}^T \bb_{ij},
         \quad \bd_{1ij} \equiv \text{first $q_2$ rows of $\bd_{0ij}$},
         \quad \bd_{2ij} \equiv \text{remaining rows of $\bd_{0ij}$}, \\
      &\D{0ij} \equiv \bQ_{ij}^T \B{ij},
         \quad \D{1ij} \equiv \text{first $q_2$ rows of $\D{0ij}$},
         \quad \D{2ij} \equiv \text{remaining rows of $\D{0ij}$},\\
      &\dD{0ij} \equiv \bQ_{ij}^T \dB{ij},
         \quad \dD{1ij} \equiv \text{first $q_2$ rows of $\dD{0ij}$},
         \quad \dD{2ij} \equiv \text{remaining rows of $\dD{0ij}$}.
 \end{align*}
\end{enumerate}
\item 
\begin{enumerate}
\item
Decompose $\stack\limits_{1 \le\, j \le n_i} (\dD{2ij}) = \bQ_i \begin{bmatrix} \bR_i \\ \bzero \end{bmatrix}$ 
such that $\bQ_i^{-1} = \bQ_i^T$ and $\bR_i$ is upper-triangular. 
\item Then obtain
\begin{align*}
      &\bc_{0i} \equiv \bQ_i^T \left\{ \stack_{1 \le\, j \le n_i} (\bd_{2ij}) \right\},
         \quad \bc_{1i} \equiv \text{first $q_1$ rows of $\bc_{0i}$},
         \quad \bc_{2i} \equiv \text{remaining rows of $\bc_{0i}$} \\
      &\bC_{0i} \equiv \bQ_i^T \left\{ \stack_{1 \le\, j \le n_i} (\D{2ij}) \right\},
         \quad \bC_{1i} \equiv \text{first $q_1$ rows of $\bC_{0i}$},
         \quad \bC_{2i} \equiv \text{remaining rows of $\bC_{0i}$}. 
\end{align*}
\end{enumerate}
\end{enumerate}
\item Decompose $\stack\limits_{1 \le i \le m} (\bC_{2i}) = \bQ 
\begin{bmatrix} \bR \\ \bzero \end{bmatrix}$ such that $\bQ^{-1} = \bQ^T$ 
and $\bR$ is upper-triangular and let
$$\TNbomega \equiv \text{first $p$ rows of } \bQ^T\left\{\stack_{1\le i\le m}(\TNbomega_{2i})\right\}.$$
\item The solutions are, for $1\le i\le m$, 
\begin{align*}
\bx_1 &= \bR^{-1}\TNbomega, \quad \Ainv{11} = \bR^{-1} \bR^{-T},\\
\bx_{2,i} &= \bR_i^{-1} (\TNbomega_{1i} - \TNbOmega_{1i}\,\bx_1), \quad
\Ainv{12,i} = -\Ainv{11} \left( \bR_i^{-1} \TNbOmega_{1i} \right)^T, \quad
\Ainv{22,i} = \bR_i^{-1}(\bR_i^{-T} - \TNbOmega_{1i}\,\Ainv{12,i}), 
\end{align*}
and, for $1\le i\le m$, $1\le j\le n_i$,
\begin{align*}
\bx_{2,ij} &= \bR_{ij}^{-1} (\bd_{1ij} - \D{1ij}\,\bx_1 - \dD{1ij}\,\bx_{2,i}), \quad
\Ainv{12,ij} = - \left\{ \bR_{ij}^{-1} (\D{1ij}\,\Ainv{11} + \dD{1ij}\,\ATinv{12,i})\right\}^T, \\
\Ainv{12,\iCOMMAj} &= - \left\{ \bR_{ij}^{-1} (\D{1ij}\,\Ainv{12,i} + \dD{1ij}\,\Ainv{22,i}) \right\}^T
\ \text{and} 
\ \Ainv{22,ij} =\bR_{ij}^{-1}(\bR_{ij}^{-T}-\D{1ij}\,\Ainv{12,ij}-\dD{1ij}\,\Ainv{12,\iCOMMAj}).
\end{align*}
\item The determinant of $\bA$ is
{\setlength\arraycolsep{3pt}
\begin{eqnarray*}
|\bA|&=&\Bigg[\big(\mbox{product of the diagonal entries of $\bR$}\big)\,\prod_{i=1}^m\,
\Bigg\{ 
\big(\mbox{product of the diagonal entries of $\bR_i$}\big)\\[1ex]
&&\qquad\qquad\qquad\qquad\qquad\qquad\qquad\ \ \ \ 
\times\prod_{j=1}^{n_i}\big(\mbox{product of the diagonal entries of $\bR_{ij}$}\big)\Bigg\}
\Bigg]^2.
\end{eqnarray*}
}
\end{enumerate}
\label{thm:fourthThm}
\end{theorem}

\noindent
Appendix \ref{sec:Th4} contains a proof of Theorem \ref{thm:fourthThm}.

\noindent
\textsl{Remarks:}
\begin{enumerate}
\item As in Theorem \ref{thm:secondThm}, Step 1(a) of Theorem \ref{thm:fourthThm} 
involves determination of $\sum_{i=1}^m n_i$ 
upper-triangular matrices $\bR_{ij}$, for $1 \le i \le m$, $1\le j\le n_i$, via QR-decomposition. 
In Step 1(b), $m$ upper-triangular matrices $\bR_{i}$, for $1 \le i \le m$, are also constructed. A final 
QR-decomposition is applied in Step 2. Each of the inversions in Step 3 can be solved rapidly via back-solving.
\item The solutions for $\bx_1,\ \Ainv{11}$ and $\bx_{2,i}, \Ainv{12,i}, \Ainv{22,i}, 1\le i\le m,$
have the same forms as in Theorem \ref{thm:fourthThm} for two-level sparse matrices. The solutions 
for $\bx_{2,ij}$, $\Ainv{12,ij}$, $\Ainv{12,\iCOMMAj}$ and $\Ainv{22,ij}$ are
suggestive of a hierarchical pattern emerging for four-level and higher level classes of the problem.
\end{enumerate}

\section{Conclusion}\label{sec:conclusion}

In this short communication we have conveyed the essence of higher level sparse matrix problems
as viewed through the prism of fitting and inference for multilevel statistical models.
Both time-honoured best linear unbiased prediction and new-fashioned mean field
variational Bayes approaches benefit from our four theorems for the two-level
and three-level situations, with details given in \citet{Nolan20}. Future extensions 
to higher level situations are aided by our results and derivations.

\section*{Acknowledgements}

We are grateful to Gioia Di Credico, Luca Maestrini, Marianne Menictas, Simon Wood
and two anonymous referees for their assistance. This research was partially supported 
by the Australian Research Council Discovery Project DP140100441.

%%%%%%%%%%%%%%%%%%%%%%%%%%%%%%%%%%%%%%%%%%%%%%%%%%%%%%%%
\appendix
%%%%%%%%%%%%%%%%%%%%%%%%%%%%%%%%%%%%%%%%%%%%%%%%%%%%%%%%%%
%
% START OF APPENDIX
%
%%%%%%%%%%%%%%%%%%%%%%%%%%%%%%%%%%%%%%%%%%%%%%%%%%%%%%%%%%

\section{Appendix: Proofs of Theorems}

%%%%%%%%%%%%%%  PROOF  OF  THEOREM  1  %%%%%%%%%%%%%%%

\subsection{Proof of Theorem 1}\label{sec:Th1}

In the case of $m=2$ the two-level sparse matrix linear system problem is
$$
\left[ 
\arraycolsep=2.2pt\def\arraystretch{1.6}
\begin{array}{c|c|c}
\setstretch{4.5}
   \A{11} & \A{12,1} & \A{12,2} \\
   \hline
   \A{12,1}^T & \A{22,1} & \bO \\
   \hline
   \A{12,2}^T & \bO & \A{22,2} \\
\end{array} 
\right]
\left[
\arraycolsep=2.2pt\def\arraystretch{1.6}
\begin{array}{c}
\setstretch{4.5}
   \bx_1 \\
   \hline
   \bx_{2,1} \\
   \hline
   \bx_{2,2} \\
\end{array} \right] =
\left[ 
\arraycolsep=2.2pt\def\arraystretch{1.6}
\begin{array}{c}
\setstretch{4.5}
   \ba_1 \\
   \hline
   \ba_{2,1} \\
   \hline
   \ba_{2,2} \\
\end{array} \right].
$$
which immediately leads to
\begin{equation*}
   \A{11} \bx_1 + \A{12,1} \bx_{2,1} + \A{12,2}\,\bx_{2,2} = \ba_1
\end{equation*}
and
\begin{equation*}
   \AT{12,i} \bx_1 + \A{22,i} \bx_{2,i} = \ba_{2,i}, \quad 1 \le i \le 2.
\end{equation*}
It is clear that the same pattern applies for general $m$, and we have
\begin{equation}
   \A{11} \bx_1 + \sum_{i=1}^m \A{12,i} \bx_{2,i} = \ba_1
\label{2levLinSyst1}
\end{equation}
and
\begin{equation}
   \AT{12,i} \bx_1 + \A{22,i} \bx_{2,i} = \ba_{2,i}, \quad 1 \le i \le m.
\label{2levLinSyst2}
\end{equation}
Conditions \eqref{2levLinSyst2} immediately imply that
\begin{equation}
   \bx_{2,i} = \A{22,i}^{-1} (\ba_{2,i} - \AT{12,i} \bx_1), \quad 1 \le i \le m.
\label{lev2x2i}
\end{equation}
Substitution of \eqref{lev2x2i} into \eqref{2levLinSyst1} then leads to the solution for
$\bx_1$ as stated in Theorem 1.

For the matrix inverse derivation, we again start with the $m = 2$ case and note that
$$
\left[ 
\arraycolsep=2.2pt\def\arraystretch{1.6}
\begin{array}{c | c | c}
\setstretch{4.5}
   \A{11} & \A{12,1} & \A{12,2} \\
   \hline
   \A{12,1}^T & \A{22,1} & \bO \\
   \hline
   \A{12,2}^T & \bO & \A{22,2} \\
\end{array}\right]
\left[ 
\arraycolsep=2.2pt\def\arraystretch{1.6}
\begin{array}{c | c | c}
\setstretch{4.5}
   \Ainv{11} & \Ainv{12,1} & \Ainv{12,2} \\
   \hline
   \ATinv{12,1} & \Ainv{22,1} & \bigtimes \\
   \hline
   \ATinv{12,2} & \bigtimes & \Ainv{22,2} \\
\end{array}\right]=
\left[ 
\arraycolsep=2.2pt\def\arraystretch{1.6}
\begin{array}{c | c | c}
\setstretch{4.5}
  \hgap\bI\hgap&\hgap\bO\hgap&\hgap\bO\hgap\\
 \hline
  \hgap\bO\hgap&\hgap\bI\hgap&\hgap\bO\hgap\\
   \hline
 \hgap\bO\hgap&\hgap\bO\hgap&\hgap\bI\hgap 
\end{array} \right].
$$
Observing the pattern from the $m=2$ case and then extending to general $m$ we obtain
the system of equations:
\begin{align}
   \A{11} \Ainv{11} + \sum_{i=1}^m \A{12,i} \ATinv{12,i} &= \bI \label{2levInv1} \\
   \AT{12,i} \Ainv{12,i} + \A{22,i} \Ainv{22,i} &= \bI, \quad 1 \le i \le m \label{2levInv2} \\
   \AT{12,i} \Ainv{11} + \A{22,i} \ATinv{12,i} &= \bO, \quad 1 \le i \le m \label{2levInv3}.
\end{align}
From \eqref{2levInv3}
\begin{equation}
   \ATinv{12,i} = - \A{22,i}^{-1} \AT{12,i} \Ainv{11}, \quad 1 \le i \le m.
\label{lev2A12i}
\end{equation}
Substitution of \eqref{lev2A12i} into \eqref{2levInv1} gives
\begin{equation*}
   \A{11} \Ainv{11} - \sum_{i=1}^m \A{12,i} \A{22,i}^{-1} \AT{12,i} \Ainv{11} = \bI
\end{equation*}
which implies that
\begin{equation*}
   \Ainv{11} = \left( \A{11} - \sum_{i=1}^m \A{12,i} \A{22,i}^{-1} \AT{12,i} \right)^{-1}.
\end{equation*}
Substitution of \eqref{lev2A12i} into \eqref{2levInv2} gives
\begin{equation*}
   - \AT{12,i} \Ainv{11} \A{12,i} \A{22,i}^{-1} + \A{22,i} \Ainv{22,i} = \bI, \quad 1 \le i \le m
\end{equation*}
implying that
$$
\Ainv{22,i} = \A{22,i}^{-1}(\bI + \AT{12,i} \Ainv{11} \A{12,i} \A{22,i}^{-1})
= \A{22,i}^{-1}(\bI - \AT{12,i}\Ainv{12,i})
, \quad 1 \le i \le m.
$$

For the $|\bA|$ result we first prove:
\vskip2mm

\noindent
\textbf{Lemma 1.} \textit{Let $\bM$ be a symmetric invertible matrix
with sub-block partitioning according to the notation}
$$\bM=\left[         
\begin{array}{cc}
\bM_{11} & \bM_{12} \\[1ex]
\bM_{12}^T & \bM_{22}
\end{array}
\right]\quad\mbox{\textit{and}}\quad
\bM^{-1}=\left[         
\begin{array}{cc}
\bM^{11} & \bM^{12} \\[1ex]
\bM^{12T} & \bM^{22}
\end{array}
\right].
$$
\textit{Then}
$$|\bM|=\left|(\bM^{11})^{-1}\right|\,\left|\bM_{22}\right|.$$
\textbf{Proof of Lemma 1.} Lemma 1 is a direct consequence of 
Theorem 13.3.8 of \citet{Harville08} concerning the determinant of a
matrix with $2\times2$ sub-block partitioning. 
\vskip5mm
\noindent
From Lemma 1 we have
$$|\bA|=\left|(\bA^{11})^{-1}\right|\left|\blockdiag_{1\le i\le m}(\bA_{22,i})\right|  
=\left|(\bA^{11})^{-1}\right|\,\prod_{i=1}^m\,|\bA_{22,i}|.
$$

%%%%%%%%%%%%%%  PROOF  OF  THEOREM  2  %%%%%%%%%%%%%%%

\subsection{Proof of Theorem 2}\label{sec:Th2}

We first note the following simplification:
$$\BT{i} \B{i} = \BT{i} \bQ_{i} \bQ_{i}^T \B{i}= \CT{0i} \C{0i}= \CT{1i} \C{1i} + \CT{2i} \C{2i},\quad 1\le i\le m,$$
where the first equality holds by the orthogonality of $\bQ_{i}$ and the second and third equalities hold
by Step 1(b) of Theorem 2. A similar sequence of steps can be used to show that
\begin{align*}
&\BT{i} \dB{i} = \CT{1i} \bR_i, \quad \dBT{i} \dB{i} = \bR_i^T \bR_i,\\
&\BT{i} \bb_i = \CT{1i} \bc_{1i} + \CT{2i} \bc_{2i} \quad \text{and} \quad \dBT{i} \bb_i = \bR_i^T \bc_{1i},\ 1\le i\le m.
\end{align*}
These simplifications allow us to represent the non-zero components of $\bA$ and the sub-vectors of $\ba$ as
\begin{equation*}
\bA_{11} = \sum_{i=1}^m (\CT{1i} \C{1i} + \CT{2i} \C{2i}), \quad 
\ba_1 = \sum_{i=1}^m (\CT{1i} \bc_{1i} + \CT{2i} \bc_{2i})
\end{equation*}
and
\begin{equation*}
   \A{12,i} = \CT{1i} \bR_i, \quad \A{22,i} = \bR_i^T \bR_i, \quad \ba_{2,i} = \bR_i^T \bc_{1i}, \quad 1 \le i \le m.
\end{equation*}

The derivation of the inverse matrix problem is as follows:
\begin{align*}
   \Ainv{11}
      &= \left( \A{11} - \sum_{i=1}^m \A{12,i} \A{22,i}^{-1} \AT{12,i} \right)^{-1} \\
      &= \left(
              \sum_{i=1}^m (\CT{1i} \C{1i} + \CT{2i} \C{2i})
              - \sum_{i=1}^m \CT{1i} \bR_i (\bR_i^T \bR_i)^{-1} \bR_i ^T \C{1i}
           \right)^{-1} \\
      &= \left( \sum_{i=1}^m \CT{2i} \C{2i} \right)^{-1} \\
      &= \left[ \left\{ \stack_{1 \le i \le m} (\C{2i}) \right\}^T 
\left\{ \stack_{1 \le i \le m} (\C{2i}) \right\} \right]^{-1} \\
      &= \left( \begin{bmatrix} \bR \\ \bO \end{bmatrix}^T \bQ^T \bQ 
\begin{bmatrix} \bR \\ \bO \end{bmatrix} \right)^{-1} \\
      &= \bR^{-1} \bR^{-T},
\end{align*}
where the fifth equality holds by Step 2 of Theorem 2, and we have used the orthogonality of $\bQ$ for the
sixth equality. The other components of $\bA^{-1}$ are found by simply substituting the above simplifications into
Theorem 1.

The derivation of the solution to the linear system is:
\begin{align*}
   \bx_1
      &= \Ainv{11} \left( \ba_1 - \sum_{i=1}^m \A{12,i} \A{22,i}^{-1} \ba_{2,i} \right) \\
      &= \bR^{-1} \bR^{-T} \left\{
              \sum_{i=1}^m (\CT{1i} \bc_{1i} + \CT{2i} \bc_{2i})
              - \sum_{i=1}^m \CT{1i} \bR_i (\bR_i^T \bR_i)^{-1} \bR_i^T \bc_{1i}
           \right\} \\
      &= \bR^{-1} \bR^{-T} \sum_{i=1}^m \CT{2i} \bc_{2i} \\
      &= \bR^{-1} \bR^{-T} \left\{ \stack_{1 \le i \le m} (\C{2i}) \right\}^T 
\left\{ \stack_{1 \le i \le m} (\bc_{2i}) \right\} \\
      &= \bR^{-1} \bR^{-T} \begin{bmatrix} \bR \\ \bO \end{bmatrix}^T \bQ^T 
\left\{ \stack_{1 \le i \le m} (\bc_{2i}) \right\} \\
      &= \bR^{-1} \bR^{-T} \bR^T \TNbomega \\
      &= \bR^{-1} \TNbomega,
\end{align*}
\noindent where we have used Step 2 of Theorem 2 for the fifth and sixth equalities. The other sub-vectors of
$\bx$ are found by simply substituting the above simplifications into Theorem 1.

For the $|\bA|$ result, the $|\bA|$ expression from Theorem 1 implies that
\begin{equation}
|\bA|=|(\Ainv{11})^{-1}|\prod_{i=1}^m\,|\A{22,i}|=|\bR^T\bR|\prod_{i=1}^m\,\big|\dBT{i} \dB{i}\big|=
|\bR^T\bR|\prod_{i=1}^m\,\big|\bR_i^T\bR_i\big|.
\label{eq:detAreR}
\end{equation}
Then
\begin{equation}
|\bR^T\bR|=|\bR^T||\bR|=(|\bR|)^2=(\mbox{product of the diagonal entries of $\bR$})^2
\label{eq:RdiagResult}
\end{equation}
where we have used the result $|\bM^T|=|\bM|$ for any square matrix $\bM$ and the fact
that the determinant of an upper-triangular matrix is the product of its diagonal
entries (Lemma 13.1.1 of \citet{Harville08}). Replacement $\bR$ with $\bR_i$ in 
(\ref{eq:RdiagResult}) and substitution into (\ref{eq:detAreR}) leads to the stated 
result for $|\bA|$.

%%%%%%%%%%%%%%  PROOF  OF  THEOREM  3  %%%%%%%%%%%%%%%

\subsection{Proof of Theorem 3}\label{sec:Th3}

In the case of $m = 2$, $n_1 = 2$ and $n_2 = 3$ the three-level sparse matrix linear system problem is
$$
\left[ 
\arraycolsep=2.2pt\def\arraystretch{1.6}
\begin{array}{c | c | c | c | c | c | c | c}
\setstretch{4.5}
   \A{11} & \A{12,1} & \A{12,11} & \A{12,12} & \A{12,2} & \A{12,21} & \A{12,22} & \A{12,23} \\
   \hline
   \A{12,1}^T & \A{22,1} & \A{12,1,1} & \A{12,1,2} & \bO & \bO & \bO & \bO \\
   \hline
   \A{12,11}^T & \A{12,1,1}^T & \A{22,11} & \bO & \bO & \bO & \bO & \bO \\
   \hline
   \A{12,12}^T & \A{12,1,2}^T & \bO & \A{22,12} & \bO & \bO & \bO & \bO \\
   \hline
   \A{12,2}^T & \bO & \bO & \bO & \A{22,2} & \A{12,2,1} & \A{12,2,2} & \A{12,2,3} \\
   \hline
   \A{12,21}^T & \bO & \bO & \bO & \A{12,2,1}^T & \A{22,21} & \bO & \bO \\
   \hline
   \A{12,22}^T & \bO & \bO & \bO & \A{12,2,2}^T & \bO & \A{22,22} & \bO \\
   \hline
   \A{12,23}^T & \bO & \bO & \bO & \A{12,2,3}^T & \bO & \bO & \A{22,23}
\end{array} \right]
\left[ 
\arraycolsep=2.2pt\def\arraystretch{1.6}
\begin{array}{c}
\setstretch{4.5}
   \bx_1 \\
   \hline
   \bx_{2,1} \\
   \hline
   \bx_{2,11} \\
   \hline
   \bx_{2,12} \\
   \hline
   \bx_{2,2} \\
   \hline
   \bx_{2,21} \\
   \hline
   \bx_{2,22} \\
   \hline
   \bx_{2,23}
\end{array} \right] =
\left[ 
\arraycolsep=2.2pt\def\arraystretch{1.6}
\begin{array}{c}
\setstretch{4.5}
   \ba_1 \\
   \hline
   \ba_{2,1} \\
   \hline
   \ba_{2,11} \\
   \hline
   \ba_{2,12} \\
   \hline
   \ba_{2,2} \\
   \hline
   \ba_{2,21} \\
   \hline
   \ba_{2,22} \\
   \hline
   \ba_{2,23}
\end{array} \right].
$$
For arbitrary values of $m$ and $\{n_i\}_{1 \le i \le m}$, we immediately 
obtain the following set of equations:
\begin{align}
   \A{11} \bx_1 + \sum_{i=1}^m \A{12,i} \bx_{2,i} + \sum_{i=1}^m \sum_{j=1}^{n_i} 
        \A{12,ij}\,\bx_{2,ij} &= \ba_1 \label{3levLinSyst1} \\
   \AT{12,i} \bx_1 + \A{22,i} \bx_{2,i} + \sum_{j=1}^{n_i} \A{12,\iCOMMAj}\,\bx_{2,ij} 
&= \ba_{2,i}, \quad 1 \le i \le m \label{3levLinSyst2} \\
   \AT{12,ij} \bx_1 + \AT{12,\iCOMMAj}\,\bx_{2,i} + \A{22,ij}\,\bx_{2,ij} 
&= \ba_{2,ij}, \quad \ 1 \le i \le m,\quad 1 \le j \le n_i. \label{3levLinSyst3}
\end{align}
Conditions \eqref{3levLinSyst3} imply that
\begin{equation}
   \bx_{2,ij} = \A{22,ij}^{-1} \left( \ba_{2,ij} - \AT{12,ij} \bx_1 - \AT{12,\iCOMMAj} \bx_{2,i} \right), 
\quad 1 \le i \le m,\quad 1 \le j \le n_i,
\end{equation}
which is the solution for $\bx_{2,ij}$ as stated in Theorem 3. Substituting this result into 
conditions \eqref{3levLinSyst2} leads to
$$
   \AT{12,i} \bx_1 + \A{22,i} \bx_{2,i}
 + \sum_{j=1}^{n_i} \A{12,\iCOMMAj} \A{22,ij}^{-1} \left( \ba_{2,ij} - \AT{12,ij} \bx_1 - \AT{12,\iCOMMAj} \bx_{2,i} \right)
      = \ba_{2,i}, \quad 1 \le i \le m.
$$
Solving this equation for $\bx_{2,i}$ and using the definitions for $\bH_{12,i}$ and $\bH_{22,i}$ 
leads to the solution for $\bx_{2,i}$ as stated in Theorem 3. Finally, 
substitution of the results for $\bx_{2,i}$ and $\bx_{2,ij}$ into condition \eqref{3levLinSyst1} gives the 
solution for $\bx_1$ stated in Theorem 3.

For the matrix inverse, we again illustrate the problem with $m=2$, $n_1 = 2$ and $n_2 = 3$:
\begin{align*}
&\left[ 
\arraycolsep=2.2pt\def\arraystretch{1.6}
\begin{array}{c | c | c | c | c | c | c | c}
\setstretch{4.5}
      \A{11} & \A{12,1} & \A{12,11} & \A{12,12} & \A{12,2} & \A{12,21} & \A{12,22} & \A{12,23} \\
      \hline
      \A{12,1}^T & \A{22,1} & \A{12,1,1} & \A{12,1,2} & \bO & \bO & \bO & \bO \\
      \hline
      \A{12,11}^T & \A{12,1,1}^T & \A{22,11} & \bO & \bO & \bO & \bO & \bO \\
      \hline
      \A{12,12}^T & \A{12,1,2}^T & \bO & \A{22,12} & \bO & \bO & \bO & \bO \\
      \hline
      \A{12,2}^T & \bO & \bO & \bO & \A{22,2} & \A{12,2,1} & \A{12,2,2} & \A{12,2,3} \\
      \hline
      \A{12,21}^T & \bO & \bO & \bO & \A{12,2,1}^T & \A{22,21} & \bO & \bO \\
      \hline
      \A{12,22}^T & \bO & \bO & \bO & \A{12,2,2}^T & \bO & \A{22,22} & \bO \\
      \hline
      \A{12,23}^T & \bO & \bO & \bO & \A{12,2,3}^T & \bO & \bO & \A{22,23}
   \end{array} \right] \\
   &\qquad\qquad \times
\left[ 
\arraycolsep=2.2pt\def\arraystretch{1.6}
\begin{array}{c | c | c | c | c | c | c | c}
\setstretch{4.5}
         \Ainv{11} & \Ainv{12,1} & \Ainv{12,11} & \Ainv{12,12} & \Ainv{12,2} & \Ainv{12,21} & \Ainv{12,22} & \Ainv{12,23} \\
         \hline
         \ATinv{12,1} & \Ainv{22,1} & \Ainv{12,1,1} & \Ainv{12,1,2} & \bigtimes & \bigtimes & \bigtimes & \bigtimes \\
         \hline
         \ATinv{12,11} & \ATinv{12,1,1} & \Ainv{22,11} & \bigtimes & \bigtimes & \bigtimes & \bigtimes & \bigtimes \\
         \hline
         \ATinv{12,12} & \ATinv{12,1,2} & \bigtimes & \Ainv{22,12} & \bigtimes & \bigtimes & \bigtimes & \bigtimes \\
         \hline
         \ATinv{12,2} & \bigtimes & \bigtimes & \bigtimes & \Ainv{22,2} & \Ainv{12,2,1} & \Ainv{12,2,2} & \Ainv{12,2,3} \\
         \hline
         \ATinv{12,21} & \bigtimes & \bigtimes & \bigtimes & \ATinv{12,2,1} & \Ainv{22,21} & \bigtimes & \bigtimes \\
         \hline
         \ATinv{12,22} & \bigtimes & \bigtimes & \bigtimes & \ATinv{12,2,2} & \bigtimes & \Ainv{22,22} & \bigtimes \\
         \hline
         \ATinv{12,23} & \bigtimes & \bigtimes & \bigtimes & \ATinv{12,2,3} & \bigtimes & \bigtimes & \Ainv{22,23}
      \end{array} \right] \\
&=\left[ 
\arraycolsep=2.2pt\def\arraystretch{1.6}
\begin{array}{c | c | c | c | c | c | c | c}
\setstretch{4.5}
         \hgap\bI\hgap&\hgap\bO\hgap&\hgap\bO\hgap&\hgap\bO\hgap&\hgap\bO\hgap&\hgap\bO\hgap&\hgap\bO\hgap&\hgap\bO\hgap\\
         \hline
         \hgap\bO\hgap&\hgap\bI\hgap&\hgap\bO\hgap&\hgap\bO\hgap&\hgap\bO\hgap&\hgap\bO\hgap&\hgap\bO\hgap&\hgap\bO\hgap\\
         \hline
         \hgap\bO\hgap&\hgap\bO\hgap&\hgap\bI\hgap&\hgap\bO\hgap&\hgap\bO\hgap&\hgap\bO\hgap&\hgap\bO\hgap&\hgap\bO\hgap\\
         \hline
         \hgap\bO\hgap&\hgap\bO\hgap&\hgap\bO\hgap&\hgap\bI\hgap&\hgap\bO\hgap&\hgap\bO\hgap&\hgap\bO\hgap&\hgap\bO\hgap\\
         \hline
         \hgap\bO\hgap&\hgap\bO\hgap&\hgap\bO\hgap&\hgap\bO\hgap&\hgap\bI\hgap&\hgap\bO\hgap&\hgap\bO\hgap&\hgap\bO\hgap\\
         \hline
         \hgap\bO\hgap&\hgap\bO\hgap&\hgap\bO\hgap&\hgap\bO\hgap&\hgap\bO\hgap&\hgap\bI\hgap&\hgap\bO\hgap&\hgap\bO\hgap\\
         \hline
         \hgap\bO\hgap&\hgap\bO\hgap&\hgap\bO\hgap&\hgap\bO\hgap&\hgap\bO\hgap&\hgap\bO\hgap&\hgap\bI\hgap&\hgap\bO\hgap\\
         \hline
         \hgap\bO\hgap&\hgap\bO\hgap&\hgap\bO\hgap&\hgap\bO\hgap&\hgap\bO\hgap&\hgap\bO\hgap&\hgap\bO\hgap&\hgap\bI\hgap
      \end{array} \right].
\end{align*}
Observing the pattern for the $m=2$, $n_1 = 2$ and $n_2 = 3$ case and extending to general 
$m$ and $\{ n_i \}_{1 \le i \le m}$, we obtain the following system of equations:
\begin{align}
   \A{11} \Ainv{11} + \sum_{i=1}^m \A{12,i} \ATinv{12,i} + \sum_{i=1}^m 
\sum_{j=1}^{n_i} \A{12,ij} \ATinv{12,ij} &= \bI \label{3levInv1} \\
   \AT{12,i} \Ainv{12,i} + \A{22,i} \Ainv{22,i} + \sum_{j=1}^{n_i} \A{12,\iCOMMAj} 
\ATinv{12,\iCOMMAj} &= \bI, \quad 1 \le i \le m, \label{3levInv2} \\
   \AT{12,ij} \Ainv{12,ij} + \AT{12,\iCOMMAj} \Ainv{12,\iCOMMAj} + \A{22,ij} 
\Ainv{22,ij} &= \bI,\quad 1 \le i \le m,\quad 1 \le j \le n_i,\label{3levInv3} \\
   \AT{12,i} \Ainv{11} + \A{22,i} \ATinv{12,i} + \sum_{j=1}^{n_i} \A{12,\iCOMMAj} 
\ATinv{12,ij} &= \bO, \quad 1 \le i \le m, \label{3levInv4} \\
   \AT{12,ij} \Ainv{11} + \AT{12,\iCOMMAj} \ATinv{12,i} + \A{22,ij} \ATinv{12,ij} 
&= \bO, \quad 1 \le i \le m,\quad 1 \le j \le n_i,\label{3levInv5} \\[1ex]
   \AT{12,ij} \Ainv{12,i} + \AT{12,\iCOMMAj} \Ainv{22,i} + \A{22,ij} \ATinv{12,\iCOMMAj} 
&= \bO,\quad 1 \le i \le m,\quad 1 \le j \le n_i.
\label{3levInv6}
\end{align}
Rearranging conditions \eqref{3levInv5} we obtain
\begin{equation*}
\Ainv{12,ij} = -\left( \Ainv{11} \A{12,ij} + \Ainv{12,i} \A{12,\iCOMMAj} \right) \A{22,ij}^{-1} 
\quad 1 \le i \le m,\quad 1 \le j \le n_i, 
\end{equation*}
from which the result for $\Ainv{12,ij}$ stated in Theorem 3 quickly follows. 
Rearranging conditions \eqref{3levInv4} we obtain
\begin{equation*}
\Ainv{12,i} = - \left( \Ainv{11} \A{12,i} + \sum_{j=1}^{n_i} \Ainv{12,ij} \AT{12,\iCOMMAj} \right) \A{22,i}^{-1}.
\end{equation*}
Substituting the results for each $\Ainv{12,ij}$ leads to
$$\Ainv{12,i} =
      - \left\{
         \Ainv{11} \A{12,i}
         - \sum_{j=1}^{n_i} \left( \Ainv{11} \A{12,ij} + \Ainv{12,i} \A{12,\iCOMMAj} \right) \A{22,ij}^{-1} \AT{12,\iCOMMAj}
      \right\} \A{22,i}^{-1},\ \ 1 \le i \le m,\ \ 1 \le j \le n_i.
$$
Solving the above set of equations for each $\Ainv{12,i}$ and using the definitions of $\bH_{12,i}$ 
and $\bH_{22,i}$ we obtain the result in Theorem 3. 
Substituting the results for each $\Ainv{12,ij}$ and $\Ainv{12,i}$ into \eqref{3levInv1}, 
using the definition of $\bH_{12,i}$ and solving for $\Ainv{11}$ leads to its stated result 
in Theorem 3.
Rearranging conditions \eqref{3levInv6} we obtain
$$
\Ainv{12,\iCOMMAj} = - (\ATinv{12,i} \A{12,ij} + \Ainv{22,i} \A{12,\iCOMMAj}) \A{22,ij}^{-1}, 
\quad 1\le i \le m,\quad 1 \le j \le n_i,
$$
from which the result for $\Ainv{12,\iCOMMAj}$ stated in Theorem 3 follows quickly. Substitution of these results into 
the corresponding  equations of condition \eqref{3levInv2} leads to
\begin{equation*}
   \AT{12,i} \Ainv{12,i} + \A{22,i} \Ainv{22,i} - \sum_{j=1}^{n_i} \A{12,\iCOMMAj} \A{22,ij}^{-1} 
(\AT{12,ij} \Ainv{12,i} + \AT{12,\iCOMMAj} \Ainv{22,i}) = \bI, \quad 1 \le i \le m.
\end{equation*}
Solving the above set of equations for each $\Ainv{22,i}$ and using the definitions of 
$\bH_{12,i}$ and $\bH_{22,i}$ in equation (7) of the main article we obtain the result in Theorem 3.
Rearrangement of \eqref{3levInv3} leads to
\begin{equation*}
   \Ainv{22,ij} = \A{22,ij}^{-1} (\bI - \AT{12,ij} \Ainv{12,ij} - \AT{12,\iCOMMAj} \Ainv{12,\iCOMMAj}), 
\quad\ 1 \le i \le m,\quad 1 \le j \le n_i.
\end{equation*}

From Lemma 1 in the proof of Theorem 1,
$$|\bA|=\left|(\bA^{11})^{-1}\right|\,\prod_{i=1}^m
\left|
\begin{array}{cc}
\A{22,i}                              & \left\{\displaystyle{\stack_{1 \le j \le n_i}}(\AT{12,\iCOMMAj})\right\}^T \\[3ex]
{\displaystyle\stack_{1 \le j \le n_i}}(\AT{12,\iCOMMAj}) &{\displaystyle\blockdiag{1\le j\le n_i}}(\A{22,ij})
\end{array}
\right|.
$$
Application of Lemma 1 and Theorem 1 to the two-level sparse matrices in the last-written expression leads to
$$\left|
\begin{array}{cc}
\A{22,i}                              & \left\{\displaystyle{\stack_{1 \le j \le n_i}}(\AT{12,\iCOMMAj})\right\}^T \\[3ex]
{\displaystyle\stack_{1 \le j \le n_i}}(\AT{12,\iCOMMAj}) &{\displaystyle\blockdiag{1\le j\le n_i}}(\A{22,ij})
\end{array}
\right|
=\left|\A{22,i}-\sum_{j=1}^{n_i}\A{12,\iCOMMAj}\A{22,ij}^{-1}\AT{12,\iCOMMAj}\right|
\prod_{j=1}^{n_i}|\A{22,ij}|
$$
and the stated result for $|\bA|$ follows immediately.

%%%%%%%%%%%%%%  PROOF  OF  THEOREM  4  %%%%%%%%%%%%%%%

\subsection{Proof of Theorem 4}\label{sec:Th4}

We first note the following re-expression:
$$
   \BT{ij} \B{ij}
      = \BT{ij} \bQ_{ij} \bQ_{ij}^T \B{ij} 
      = \begin{bmatrix} \D{1ij} \\ \D{2ij} \end{bmatrix}^T \begin{bmatrix} \D{1ij} \\ \D{2ij} \end{bmatrix} 
      = \DT{1ij} \D{1ij} + \DT{2ij} \D{2ij},\ 1\le i\le m,\quad 1\le j\le n_i,
$$
where the first equality holds by orthogonality of $\bQ_{ij}$ and the second equality 
holds by Step 1(a)ii of Theorem 4. A similar sequence of steps can be used to show that
\begin{align*}
   &\dBT{ij} \dB{ij} = \dDT{1ij} \dD{1ij} + \dDT{2ij} \dD{2ij}, \quad \BT{ij} \dB{ij} 
= \DT{1ij} \dD{1ij} + \DT{2ij} \dD{2ij},  \\
   &\BT{ij} \bb_{ij} = \DT{1ij} \bd_{1ij} + \DT{2ij} \bd_{2ij} \quad 
\text{and} \quad \dBT{ij} \bb_{ij} = \dDT{1ij} \bd_{1ij} + \dDT{2ij} \bd_{2ij},\quad1\le i\le m,\quad1\le j\le n_i. 
\end{align*}
We also have
$$\ddBT{ij} \ddB{ij}= \ddBT{ij} \bQ_{ij} \bQ_{ij}^T \ddB{ij} 
= \begin{bmatrix} \bR_{ij} \\ \bzero \end{bmatrix}^T \begin{bmatrix} \bR_{ij} \\ \bzero \end{bmatrix} \\
= \bR_{ij}^T \bR_{ij},
$$
where the first equality holds by orthogonality of $\bQ_{ij}$ and the second equality holds by 
Step 1(a)i of Theorem 4. A similar sequence of steps can be used to show that
\begin{align*}
   \dBT{ij} \ddB{ij} = \dDT{1ij} \bR_{ij}, \quad \BT{ij} \ddB{ij} = \DT{1ij} \bR_{ij} 
\quad \text{and} \quad \ddBT{ij} \bb_{ij} = \bR_{ij}^T \bd_{1ij}.
\end{align*}
The above simplifications allow us to represent the non-zero components of $\bA$ and the sub-vectors of $\ba$ as
\begin{align*}
   &\A{11} = \sum_{i=1}^m \sum_{j=1}^{n_i} (\DT{1ij} \D{1ij} + \DT{2ij} \D{2ij}), \quad
      \ba_1 = \sum_{i=1}^m \sum_{j=1}^{n_i} (\DT{1ij} \bd_{1ij} + \DT{2ij} \bd_{2ij}), \\
   &\A{22,i} = \sum_{j=1}^{n_i} (\dDT{1ij} \dD{1ij} + \dDT{2ij} \dD{2ij}), \quad
      \A{12,i} = \sum_{j=1}^{n_i} (\DT{1ij} \dD{1ij} + \DT{2ij} \dD{2ij}) \\
   &\ba_{2,i} = \sum_{j=1}^{n_i} (\dDT{1ij} \bd_{1ij} + \dDT{2ij} \bd_{2ij}), \quad 1 \le i \le m,
\end{align*}
and, for $1\le i \le m,\ 1 \le j \le n_i$,
$$\A{22,ij} = \bR_{ij}^T \bR_{ij}, \quad \A{12,\iCOMMAj} = \dDT{1ij} \bR_{ij},\quad
\A{12,ij} = \DT{1ij} \bR_{ij}, \quad \ba_{2,ij} = \bR_{ij}^T \bd_{1ij}.
$$
Furthermore, each $\bH_{22,i}$ matrix takes the form
\begin{align*}
   \bH_{22,i}
      &= \A{22,i} - \sum_{j=1}^{n_i} \A{12,i \ j} \A{22,ij}^{-1} \AT{12,\iCOMMAj} \\
      &= \sum_{j=1}^{n_i} (\dDT{1ij} \dD{1ij} + \dDT{2ij} \dD{2ij}) 
- \sum_{j=1}^{n_i} \dDT{1ij} \bR_{ij} (\bR_{ij}^T \bR_{ij})^{-1} \bR_{ij}^T \dD{1ij} \\
      &= \sum_{j=1}^{n_i} \dDT{2ij} \dD{2ij} \\
      &= \left\{ \stack_{1 \le j \le n_i} \dD{2ij} \right\}^T \left\{ \stack_{1 \le j \le n_i} \dD{2ij} \right\} \\
      &= \left( \bQ_i \begin{bmatrix} \bR_i \\ \bzero \end{bmatrix} \right)^T \bQ_i 
\begin{bmatrix} \bR_i \\ \bzero \end{bmatrix}\\
      &= \begin{bmatrix} \bR_i \\ \bzero \end{bmatrix}^T \begin{bmatrix} \bR_i \\ \bzero \end{bmatrix} \\
      &= \bR_i^T \bR_i, \quad 1 \le i \le m,
\end{align*}
where the fifth equality holds by Step 1(b)ii of Theorem 4. Similarly,
$$\bH_{12,i} = \TNbOmega_{1i}^T \bR_i\quad\mbox{and}\quad 
\bH_{12,i} \bH_{22,i}^{-1} \bH_{12,i}^T = \TNbOmega_{1i}^T \TNbOmega_{1i}, \quad 1 \le i \le m.$$

We are now in a position to derive the solutions to the inverse matrix problem:
\begin{align*}
   \Ainv{11}
      &= \left\{
               \sum_{i=1}^m \sum_{j=1}^{n_i} (\DT{1ij} \D{1ij} + \DT{2ij} \D{2ij})
               - \sum_{i=i}^m \TNbOmega_{1i}^T \TNbOmega_{1i}
               - \sum_{i=1}^m\sum_{j=1}^{n_i} \DT{1ij} \bR_{ij} (\bR_{ij}^T \bR_{ij})^{-1} \bR_{ij}^T \D{1ij}
            \right\}^{-1} \\
      &= \left(
               \sum_{i=1}^m \sum_{j=1}^{n_i} \DT{2ij} \D{2ij}
               - \sum_{i=i}^m \TNbOmega_{1i}^T \TNbOmega_{1i}
            \right)^{-1} \\
      &= \left(
               \sum_{i=1}^m \left[
                  \left\{ \stack_{1 \le j \le n_i} (\D{2ij}) \right\}^T 
                 \left\{ \stack_{1 \le j \le n_i} (\D{2ij}) \right\}
                  - \TNbOmega_{1i}^T \TNbOmega_{1i}
               \right]
            \right)^{-1} \\
      &= \left(
               \sum_{i=1}^m \left[
                  \left\{ \stack_{1 \le j \le n_i} (\D{2ij}) \right\}^T \bQ_i \bQ_i^T 
             \left\{ \stack_{1 \le j \le n_i} (\D{2ij}) \right\}
                  - \TNbOmega_{1i}^T \TNbOmega_{1i}
               \right]
            \right)^{-1} \\
      &= \left(
               \sum_{i=1}^m \left[
                  \begin{bmatrix} \TNbOmega_{1i} \\ \TNbOmega_{2i} \end{bmatrix}^T 
            \begin{bmatrix} \TNbOmega_{1i} \\ \TNbOmega_{2i} \end{bmatrix}
                  - \TNbOmega_{1i}^T \TNbOmega_{1i}
               \right]
            \right)^{-1} \\
      &= \left(
               \sum_{i=1}^m \TNbOmega_{2i}^T \TNbOmega_{2i}
            \right)^{-1} \\
      &= \left[
               \left\{ \stack_{1 \le i \le m} (\TNbOmega_{2i}) \right\}^T 
\left\{ \stack_{1 \le i \le m} (\TNbOmega_{2i}) \right\}
            \right]^{-1} \\
      &= \left\{
               \left( \bQ \begin{bmatrix} \bR \\ \bzero \end{bmatrix} \right)^T \bQ 
\begin{bmatrix} \bR \\ \bzero \end{bmatrix}
            \right\}^{-1} \\
      &= \left(
               \begin{bmatrix} \bR \\ \bzero \end{bmatrix}^T \begin{bmatrix} \bR \\ \bzero \end{bmatrix}
            \right)^{-1} \\
      &= \bR^{-1} \bR^{-T},
\end{align*}
where the fourth equality holds by orthogonality of each $\bQ_i$, the fifth equality holds by Step 1(b)ii, 
the eighth equality holds by Step 2 and the ninth equality holds by orthogonality of $\bQ$. The other 
components of $\bA^{-1}$ are found by simply substituting the above simplifications into Theorem 3
and using a similar sequence of steps.

For the linear system solution we have
\begin{align*}
 \bh_{2,i}
&= \sum_{j=1}^{n_i} (\dDT{1ij} \bd_{1ij} + \dDT{2ij} \bd_{2ij}) - \sum_{j=1}^{n_i} \dDT{1ij} 
\bR_{ij} (\bR_{ij}^T \bR_{ij})^{-1} \bR_{ij}^T \bd_{1ij} 
= \sum_{j=1}^{n_i} \dDT{2ij} \bd_{2ij}\\
&=\left\{ \stack_{1 \le j \le n_i} (\dD{2ij}) \right\}^T \left\{ \stack_{1 \le j \le n_i} (\bd_{2ij}) \right\}
=\begin{bmatrix} \bR_i \\ \bzero \end{bmatrix}^T\bQ_i^T\left\{ \stack_{1 \le j \le n_i} (\bd_{2ij}) \right\} \\
&=\begin{bmatrix} \bR_i \\ \bzero \end{bmatrix}^T\begin{bmatrix} \TNbomega_{1i} \\ \TNbomega_{2i} \end{bmatrix}
=\bR_i^T\TNbomega_{1i}, \quad 1 \le i \le m,
\end{align*}
where the fourth and fifth equalities hold by Step 1(b)i and Step 1(b)ii, respectively. Also,
$$\A{12,ij} \A{22,ij}^{-1} \ba_{2,ij}
= \DT{1ij} \bR_{ij} (\bR_{ij}^T \bR_{ij})^{-1} \bR_{ij}^T \bd_{1ij} 
= \DT{1ij} \bd_{1ij}, \quad 1 \le i \le m,\quad 1 \le j \le n_i. 
$$
The final task is to derive the expressions in Theorem 4 for the components
of $\bx$. For the first block we have:
\begin{align*}
   \bx_1
      &= \Ainv{11} \left\{
               \sum_{i=1}^m \sum_{j=1}^{n_i} (\DT{1ij}\bd_{1ij} + \DT{2ij}\bd_{2ij})
               - \sum_{i=1}^m \TNbOmega_{1i}^T \bR_i^{-T} \bR_i^{T}\TNbomega_{1i}
               - \sum_{i=1}^m \sum_{j=1}^{n_i} \DT{1ij}\bd_{1ij}
            \right\} \\
      &= \Ainv{11} \left\{
               \sum_{i=1}^m \left(
                  \sum_{j=1}^{n_i} \DT{2ij}\bd_{2ij}
                  -  \TNbOmega_{1i}^T \TNbomega_{1i}
               \right)
            \right\} \\
      &= \Ainv{11} \left(
               \sum_{i=1}^m \left[
                  \left\{ \stack_{1 \le j \le n_i} (\D{2ij}) \right\}^T \left\{ \stack_{1 \le j \le n_i} (\bd_{2ij}) \right\}
                  -  \TNbOmega_{1i}^T \TNbomega_{1i}\right]\right) \\
      \begin{split}
         &= \Ainv{11} \left(
                  \sum_{i=1}^m \left[
                     \left\{ \stack_{1 \le j \le n_i} (\D{2ij}) \right\}^T \bQ_i \bQ_i^T 
\left\{ \stack_{1 \le j \le n_i} (\bd_{2ij}) \right\} - \TNbOmega_{1i}^T \TNbomega_{1i}
                  \right]
               \right)
      \end{split}\\
      &= \Ainv{11} \left(
               \sum_{i=1}^m \left[
                  \begin{bmatrix} \TNbOmega_{1i} \\ \TNbOmega_{2i} \end{bmatrix}^T 
\begin{bmatrix} \TNbomega_{1i} \\ \TNbomega_{2i} \end{bmatrix}
                  - \TNbOmega_{1i}^T \TNbomega_{1i}
               \right]
            \right) \\
      &= \Ainv{11} \left\{
               \sum_{i=1}^m \left(
                  \TNbOmega_{1i}^T \TNbomega_{1i} + \TNbOmega_{2i}^T \TNbomega_{2i}
                  - \TNbOmega_{1i}^T\TNbomega_{1i}
               \right)
            \right\} \\
      &= \Ainv{11} \sum_{i=1}^m \TNbOmega_{2i}^T \TNbomega_{2i} \\
      &= \Ainv{11} \left\{ \stack_{1 \le i \le m} (\TNbOmega_{2i}) \right\}^T 
\left\{ \stack_{1 \le i \le m} (\TNbomega_{2i}) \right\} \\
      &= \Ainv{11} \left( \bQ \begin{bmatrix} \bR \\ \bzero \end{bmatrix} \right)^T 
\left\{ \stack_{1 \le i \le m} (\TNbomega_{2i}) \right\} \\
      &= \bR^{-1} \bR^{-T} \begin{bmatrix} \bR \\ \bzero \end{bmatrix}^T \bQ^T 
\left\{ \stack_{1 \le i \le m} (\TNbomega_{2i}) \right\} \\
      &= \bR^{-1} \bR^{-T} \bR^T \TNbomega \\
      &= \bR^{-1} \TNbomega,
\end{align*}
where the fourth equality holds by the orthogonality of $\bQ_{i}$, the fifth equality holds by Step 1(b)
and the ninth and eleventh equalities hold by Step 2 of Theorem 4. We get the 
other sub-vectors of $\bx$ by a similar sequence of computations.

From Theorem 3 we have
$$
|\bA|=
\left|\left(\Ainv{11}\right)^{-1}\right|
\,\prod_{i=1}^m\left(\left|\bA_{22,i}-\sum_{j=1}^{n_i}\bA_{12,\iCOMMAj}\,
\A{22,ij}^{-1}\,\bA_{12,\iCOMMAj}^T\right|
\prod_{j=1}^{n_i}\big|\A{22,ij}\big|\right).
$$
From part 3. of Theorem 4 we have $\Ainv{11}=\bR^{-1}\bR^{-T}$ so steps similar to those used
for simplification of $|\bA|$ in the proof of Theorem 2 lead to 
$$\left|\left(\Ainv{11}\right)^{-1}\right|=\big(\mbox{product of the diagonal entries of $\bR$}\big)^2.$$
From the above arguments involving the $\bH_{22,i}$ matrices, 
{\setlength\arraycolsep{3pt}
\begin{eqnarray*}
\left|\bA_{22,i}-\sum_{j=1}^{n_i}\bA_{12,\iCOMMAj}\,
\A{22,ij}^{-1}\,\bA_{12,\iCOMMAj}^T\right|&=&|\bH_{22,i}|=|\bR_i^T\bR_i|
=\big(|\bR_i|\big)^2\\[1ex]
&=&\big(\mbox{product of the diagonal entries of $\bR_i$}\big)^2
\end{eqnarray*}
}
for $1\le i\le m$. Lastly, 
$$\big|\A{22,ij}\big|=|\ddBT{ij} \ddB{ij}|=|\bR_{ij}^T\bR_{ij}|=\big(|\bR_{ij}|\big)^2
=\big(\mbox{product of the diagonal entries of $\bR_{ij}$}\big)^2
$$
for $1\le i\le m,\ 1\le j\le n_i$. The stated result for $|\bA|$ immediately follows.

% In using BibteX, use wb_stat.bst
%\bibliographystyle{wb_stat}
%\bibliography{bibfile}

\end{document}